\documentstyle{article}

\input amssym.def

\newsymbol \ltimes 226E %% semi direct product signs
\newsymbol \rtimes 226F

\def\mapupalt#1{{\Big\uparrow\llap{$\vcenter {\hbox {$#1$}}$~~}}}
\def\mapdown#1{{\Big\downarrow\rlap{$\vcenter {\hbox {$#1$}}$}}}
\def\mapright#1{{\smash{\mathop{\longrightarrow}\limits^{#1}}}}
\def\mapne#1{{\nearrow\llap{$\vcenter {\hbox {$#1$}}$~~~~}}}
\def\mapse#1{{\searrow\llap{$\vcenter {\hbox {$#1$}}$~~~~}}}

\def\sqr#1#2{{\vcenter{\vbox{\hrule height.#2pt
        \hbox{\vrule width.#2pt height#1pt \kern#1pt
           \vrule width.#2pt}
        \hrule height.#2pt}}}}
\def\square{\mathchoice\sqr56\sqr56\sqr{2.1}3\sqr{1.5}3}

\def\proof{\noindent{\bf Proof.~}}
\def\endproof{\quad $\square$}

\newtheorem{thm}{Theorem}[section]
\newtheorem{lem}[thm]{Lemma}
\newtheorem{prop}[thm]{Proposition}
\newtheorem{cor}[thm]{Corollary}
\newtheorem{defn}[thm]{Definition}
\newtheorem{pagph}[thm]{($\!$}
\newtheorem{ex}[thm]{Example}
\newtheorem{rk}[thm]{Remark}
\newtheorem{conj}[thm]{Conjecture}

\newenvironment{theorem}{\begin{thm}\sl}{\end{thm}\rm}
\newenvironment{lemma}{\begin{lem}\sl}{\end{lem}\rm}
\newenvironment{proposition}{\begin{prop}\sl}{\end{prop}\rm}
\newenvironment{corollary}{\begin{cor}\sl}{\end{cor}\rm}
\newenvironment{definition}{\begin{defn}\rm}{\end{defn}}
\newenvironment{para}{\begin{pagph}$\!$\bf)\rm}{\end{pagph}}

\def\Z{{\Bbb Z}}
\def\Q{{\Bbb Q}}

\def\C{{\Bbb C}}
\def\P{{\Bbb P}}

\def\A{{\cal A}}
\def\M{{\cal M}}

\def\calB{{\cal B}}
\def\calD{{\cal D}}
\def\calE{{\cal E}}
\def\calG{{\cal G}}
\def\calH{{\cal H}}
\def\calJ{{\cal J}}
\def\calK{{\cal K}}
\def\calL{{\cal L}}
\def\calO{{\cal O}}
\def\calP{{\cal P}}
\def\calQ{{\cal Q}}
\def\calT{{\cal T}}
\def\calU{{\cal U}}
\def\calV{{\cal V}}
\def\calZ{{\cal Z}}

\def\rhotilde{\tilde{\rho}}
\def\rhohat{\hat{\rho}}
\def\Gammahat{\hat{\Gamma}}

\def\Jhat{\widehat{J}}

\def\u{{\goth u}}
\def\p{{\goth p}}

\def\Hom{{\rm Hom}}

\def\invlim{\lim_\leftarrow}
\def\im{{\rm im\, }}
\def\Aut{{\rm Aut\,}}

\def\spec{{\rm Spec\,}}
\def\Jac{{\rm Jac\,}}
\def\Pic{{\rm Pic}}
\def\Ext{{\rm Ext}}

\def\im{{\rm im\,}}

\def\semi{\ltimes} % change this to semidirect product.

\def\s{\overline{s}}

\def\comp{~\widehat{\!}{\;}}
\def\comptensor{\widehat\otimes}

\begin{document}

\begin{center}

{\bf \LARGE Completions of mapping class groups \\
and the cycle $C - C^-$}\\
\bigskip

July,  1992 \\
\bigskip

Richard~M.~Hain\footnote{Supported in part by grants from the
National
Science Foundation.}\\
\medskip

Department of Mathematics\\
Duke University\\
Durham, NC 27706
\end{center}

\section{Introduction}\label{intro}

The classical Malcev (or unipotent completion) of an abstract group
$\pi$ is a prounipotent group $\calP$ defined over $\Q$ together with
a homomorphism $\phi : \pi \to \calP$. It is characterized by the
property that if $\psi : \pi \to \calU$ is a homomorphism of $\pi$
into a prounipotent group, then there is a unique homomorphism
$\Psi: \calP \to \calU$ of prounipotent groups such that
$\psi = \Psi \phi$.
$$
\matrix{
\pi &\stackrel{\phi}{\to} &\calP \cr
& \mapse{\psi} &\mapdown{\Psi} \cr
&& \calU \cr }
$$

When $H_1(\pi;\Q) = 0$, the unipotent completion is trivial, and it
is
therefore a useless tool for studying mapping class groups. Deligne
has suggested a notion of {\it relative Malcev completion}:  Suppose
$\Gamma$ is an abstract group and that $\rho : \Gamma \to S$ is a
homomorphism of $\Gamma$ into a semisimple linear algebraic
group defined over $\Q$. Suppose that $\rho$ has
Zariski dense image. The {\it completion of $\Gamma$ relative to
$\rho$}
is a proalgebraic group $\calG$ over $\Q$, which is an extension of
$S$ by a prounipotent group $\calU$, and a homomorphism
$\tilde{\rho} : \Gamma \to \calG$ which lifts $\rho$.  When $S$ is
the trivial group, it reduces to the unipotent completion.  The
relative
completion is characterized by a universal mapping property which
generalizes the one in the unipotent case (see (\ref{ump})).

Denote the mapping class group associated to a surface of genus $g$
with $r$ boundary components and $n$ ordered marked distinct
points by $\Gamma_{g,r}^n$.  The mapping class group
$\Gamma_{g,r}^n$ has a natural representation
$\rho : \Gamma_{g,r}^n \to Sp_g(\Z)$ obtained from the action of
$\Gamma_{g,r}^n$ on the first homology group of the underlying
compact Riemann surface.  Its kernel is, by definition, the Torelli
group $T_{g,r}^n$.  One can therefore form the  completion of
$\Gamma_{g,r}^n$ relative to $\rho$.  It is a proalgebraic group
$\calG_{g,r}^n$ which is an extension
$$
1 \to \calU_{g,r}^n \to \calG_{g,r}^n \to Sp_g \to 1
$$
of $Sp_g$ by a prounipotent group.   The homomorphism
$$
\tilde{\rho}: \Gamma_{g,r}^n \to \calG_{g,r}^n
$$
induces a homomorphism $T_{g,r}^n \to \calU_{g,r}^n$. This, in turn,
induces a homomorphism $\calT_{g,r}^n \to \calU_{g,r}^n$ from the
unipotent completion of the Torelli group into $\calU_{g,r}^n$. Our
main result is:
\medskip

\noindent{\bf Theorem.}{\sl When $g\ge 3$, the natural homomorphism
$\calT_{g,r}^n \to \calU_{g,r}^n$ is surjective with nontrivial
kernel
which is contained in the center of $\calT_{g,r}^n$
and which is isomorphic to $\Q$ whenever $g \ge 8$.}
\medskip

We prove this by relating the central extension above to the line
bundle over the moduli space of genus three curves associated to the
archimedean height of the algebraic cycle $C - C^-$ in the jacobian
of
a curve $C$ of genus 3.\footnote{When the genus $g$ of $C$ is $\ge
3$,
one can relate this
central extension to the height pairing between the
cycles $C^{(a)} - {C^{(a)}}^-$ and $C^{(b)} - {C^{(b)}}^-$ in $\Jac
C$,
where $a+b = g-1$ and $C^{(r)}$ denotes the $r$th symmetric power of
$C$.
We chose not to do this in order to keep the Hodge theory
straightforward.}
This theorem is related to, and complements, the work of Morita
\cite{morita}. The constant 8 in the theorem can surely be improved,
possibly to 3.\footnote{The optimal constant is the smallest integer
$d$
such that $H^2(Sp_g(\Z),A)$ vanishes for all rational representations
$A$ of
$Sp_g(\Q)$ whenever $g\ge d$.}

One reason for introducing relative completions of fundamental
groups of varieties, and of mapping class groups in particular, is
that
their coordinate rings are, under suitable conditions, direct limits
of
variations of mixed Hodge structure over the variety.  This result
and
some of its applications to the action of the mapping class group of
a surface
$S$ on the lower central series of $\pi_1(S,\ast)$ will be presented
elsewhere.

Part of a general theory of relative completions is worked out in
Section \ref{basic}.  Many of the results of that section were worked
out independently and contemporaneously by Eduard Looijenga.
I would  like to thank him for his correspondence.
I would also like to thank P.~Deligne for his correspondence, and the
Mathematics Department of the University of Utrecht for its
hospitality and support during a visit in May, 1992 when this paper
was written.
\medskip

\noindent{\bf Conventions:} The group $Sp_g(R)$ will denote the
group of automorphisms of a free $R$ module of dimension $2g$
which preserve a unimodular skew symmetric bilinear form. In
short, elements of $Sp_g(R)$ are $2g \times 2g$ matrices.

\section{Relative  completion}
\label{rel_malcev}

Fix a field $F$ of characteristic zero.
Suppose that $\pi$ is a abstract group and that $\rho : \pi \to S$ is
a
representation of $\pi$ into a linear algebraic group $S$ defined
over
$F$.  Assume that the image of $\rho$ is Zariski
dense.  In this section we define the {\it  completion of $\pi$
relative to $\rho$}.  When $S$ is the trivial group, this reduces to
the
Malcev completion (a.k.a.\ unipotent completion) which is
defined, for example, in \cite{quillen}, \cite{chen:malcev} and
\cite{sullivan}.  The idea of relative  completion is due to Deligne
\cite{deligne:variations} and is a refinement of the idea of the
``algebraic
hull'' of a group introduced by Hochschild and Mostow
\cite[p.~1140]{hoch-mostow}.

To construct the relative completion of $\pi$ with respect to $\rho$,
consider all  commutative diagrams of the form
$$
\matrix{
1 & \to & U & \to & E & \to & S & \to &1 \cr
   &       &     &      & \rhotilde\uparrow & \nearrow\rho &&&\cr
&&&& \pi &&&&\cr}
$$
where $E$ is a linear algebraic group over $F$, $U$ a unipotent
subgroup of $E$, and where $\rhotilde$ is a lift of
$\rho$ to $E$ whose image is Zariski dense.  All
morphisms in the top row are algebraic group homomorphisms.  One can
define morphisms of such diagrams in the obvious way.

\begin{proposition} \label{invsys} The set of such diagrams forms an
inverse
system.
\end{proposition}

\proof  If
$$
\matrix{
1 & \to & U_\alpha & \to & E_\alpha & \to & S & \to &1 \cr
&&&&&&&&\cr
1 & \to & U_\beta & \to & E_\beta & \to & S & \to &1 \cr}
$$
are two extensions of $S$ by a unipotent algebraic group, then one
can form
the  fibered product
$$
\matrix{
E & \to& E_\alpha \cr
\downarrow &&\downarrow \cr
E_\beta & \to & S .\cr}
$$
The natural homomorphism $E \to S$ is surjective with kernel
the unipotent group $U_\alpha \times U_\beta$.

Now  suppose that $\rho_\alpha : \pi \to E_\alpha$ and $\rho_\beta :
\pi \to
E_\beta$  are lifts of $\rho : \pi \to S$ to $E_\alpha$ and
$E_\beta$,
respectively, both with Zariski dense image.  They induce a
homomorphism
$\rho_{\alpha\beta} : \pi \to E$ which lifts both $\rho_\alpha$ and
$\rho_\beta$.  Denote the Zariski closure of the image of
$\rho_{\alpha\beta}$ in $E$ by $E_{\alpha\beta}$. Then
$E_{\alpha\beta}$
is a linear algebraic group and the kernel of the natural
homomorphism
$E_{\alpha\beta} \to S$ is unipotent as it is a subgroup of $U_\alpha
\times U_\beta$.  The natural map $E_{\alpha\beta}\to S$ is
surjective as
the image of $\rho$ is Zariski dense in $S$. \endproof

\begin{definition}  The {\it completion $\calP_F$ of $\pi$ (over $F$)
relative to $\rho : \pi \to S$} is defined to be
the proalgebraic group
$$
\calP_F = \invlim E,
$$
where the inverse limit is taken over all commutative diagrams
$$
\matrix{
1 & \to & U & \to & E & \to & S & \to &1 \cr
   &       &     &      & \mapupalt{\rhotilde} & \mapne{\rho} &&&\cr
&&&& \pi &&&&\cr}
$$
whose top row is an extension of $S$ by a unipotent group in the
category of linear algebraic groups over $F$, and where
$\rhotilde$ has dense image.  The homomorphisms $\rhotilde : \pi
\to E$ induce a canonical homomorphism $\pi \to \calP_F$.
\end{definition}

Often we will simply say that $\pi \to \calP_F$ is the {\it $S$-%
completion\/} of $\pi$.  The coordinate ring of $\calP_F$ is the
direct limit
of
the coordinate rings of the groups $E$. It will be denoted
$\calO(\calP_F)$.
It
is a commutative Hopf algebra with antipode.

There is a natural surjection $\calP_F \to S$ whose kernel is a
prounipotent
group.  When $S$ is the trivial group, we obtain the classical Malcev
completion.  The $S$-Malcev completion is characterized by the
following easily verified  universal mapping property.

\begin{proposition} \label{ump} Suppose that $\calE$ is a linear
proalgebraic
group defined over $F$, and that $\calE \to S$ is a homomorphism of
proalgebraic groups with prounipotent kernel.  If $\varphi : \pi \to
\calE$ is
a group
homomorphism, then there is a unique homomorphism $\tau :
\calP_F \to
\calE$ of pro-algebraic groups over $F$ such that the
diagram
$$
\matrix{
&&\calP_F&& \cr
&\rhohat\nearrow& &\searrow& \cr
\pi & & \mapdown{\tau} && S \cr
&\varphi \searrow && \nearrow& \cr
&& \calE&& \cr}
$$
commutes. \endproof
\end{proposition}

Suppose that $G$ is a (pro)algebraic group over the field $F$.
Suppose that $k$ is a field extension of $F$. We shall denote $G$,
viewed as an
algebraic group over $k$ by extension of scalars, by $G(k)$. The
following assertion follows directly from the universal mapping
property.

\begin{corollary} If $k$ is a field extension of the field $F$, then
there is a
natural homomorphism $\calP_{k} \to \calP_F(k)$.
\endproof
\end{corollary}

We will show in the next two sections that, with some extra
hypotheses,
this homomorphism is always an isomorphism.

\section{A construction of the Malcev completion}
\label{classical}

There is an explicit algebraic construction of the Malcev completion
which is due to Quillen \cite{quillen}. We will use it to show that
the Malcev completion of a group over $F$ is isomorphic to
the $F$ points of its Malcev completion over $\Q$.

Denote the group algebra of a
group $\pi$ over a commutative ring $R$ by $R\pi$. The
{\it augmentation\/} is
the homomorphism $\epsilon : R\pi \to R$ defined by taking each
$\gamma\in \pi$ to 1.  The kernel of the augmentation is called the
{\it
augmentation ideal\/} and will be denoted by $J_R$, or simply $J$
when there is no chance of confusion.  With the coproduct
$\Delta : R\pi \to R\pi \otimes R\pi$ defined by
$\Delta (\gamma) = \gamma \otimes \gamma$, for all $\gamma \in \pi$,
$R\pi$ has
the structure
of a cocommutative Hopf algebra.

The powers of the augmentation ideal define a topology on $R\pi$
that is called the {\it $J$-adic topology}.  The {\it $J$-adic
completion\/} of the group ring is the $R$ module
$$
R\pi\comp := \invlim R\pi/J^l.
$$
The completion of $J$ will be denoted by $\Jhat$. Since the coproduct
is continuous, it induces a coproduct
$$
\Delta : R\pi\comp \to R\pi \comptensor R\pi,
$$
where $\comptensor$ denotes the completed tensor product.  This
gives $R\pi\comp$ the structure of  a complete Hopf algebra.

The proof of the following proposition is straightforward.

\begin{proposition} \label{first-quot}  If $\pi$ is a group and $R$ a
ring, then the function
$
\pi \to J_R/J_R^2
$
defined by taking $\gamma \in \pi$ to the coset of $\gamma-1$,
induces an
$R$-module isomorphism
$$
H_1(\pi,R) \approx J_R / J_R^2.\quad \square
$$
\end{proposition}

Now let $R$ be a field $F$ of characteristic zero. The
logarithm and exponential maps are mutually inverse
homeomorphisms
$$
\log : 1 + \Jhat_F \to \Jhat_F \hbox{ and } \exp : \Jhat_F \to
1+ \Jhat_F.
$$

The set of {\it primitive elements\/} $\p$ of $F\pi\comp$ is defined
by
$$
\p = \{ X \in \Jhat_F : \Delta(X) = X \otimes 1 + 1 \otimes X  \}.
$$
With the bracket $[X,Y] = XY - YX$,  it has the structure of a Lie
algebra.  The topology of $F\pi\comp$ induces a topology on $\p$,
giving it the structure of a complete topological Lie algebra.

The set of {\it group-like\/} elements $\calP$ of $F\pi\comp$ is
defined by
$$
\calP = \{ X \in F\pi\hat {\ } : \Delta(X) = X \otimes X
\hbox{ and } \epsilon(X) = 1 \}.
$$
It is a subgroup of the group of units of $F\pi\comp$.
The logarithm and exponential maps restrict to  mutually inverse
homeomorphisms
$$
\log : \calP \to \p \hbox{ and } \exp : \p \to \calP.
$$

The filtration of $F\pi\comp$ by the powers of $J$ induces
filtrations
of $\calP$ and $\p$:  Set
$$
\p^l = \p \cap \Jhat^l \hbox{ and } \calP^l = \calP \cap (1 +
\Jhat^l).
$$
These satisfy
$$
\p = \p^1 \supseteq \p^2 \supseteq \p^3 \cdots
$$
and
$$
\calP = \calP^1 \supseteq \calP^2 \supseteq \calP^3 \cdots .
$$
These filtrations define topologies on $\p$ and $\calP$. Both are
separated and
complete.  For each $l$ set
$$
\p_l = \p / \p^{l+1} \hbox{ and } \calP_l = \calP / \calP^{l+1}.
$$
It follows easily from (\ref{first-quot}) that if $H_1(\pi,\Q)$ is
finite dimensional (e.g. $\pi$ is finitely generated), then each
$\calP_l$ is a linear algebraic group.

Since the logarithm and exponential maps induce isomorphisms
between $1+\Jhat/\Jhat^{l+1}$ and $\Jhat/\Jhat^{l+1}$, and since
$$
\p_l \subseteq \Jhat/\Jhat^{l+1} \hbox{ and }
\calP_l \subseteq 1 + \Jhat/\Jhat^{l+1},
$$
it follows that the logarithm and exponential maps induce polynomial
bijections between $\p_l $ and $\calP_l$.  Consequently, when
$H_1(\pi,\Q)$ is finite dimensional, each of the groups
$\calP_l$ is a unipotent algebraic group over $F$ with Lie algebra
$\p_l$.
It follows that if $H_1(\pi,F)$ is finite dimensional, then $\calP$
is a
prounipotent group over $F$.

The composition of the canonical inclusion of $\pi$ into $F\pi$
followed by the completion map $F\pi \to F\pi\comp$ yields  a
canonical map $ \pi \to F\pi\comp$. Since the image of this map is
contained in $\calP$, there is a canonical homomorphism $\pi \to
\calP$.
The composition of the natural homomorphism $\pi \to \calP$ with
the quotient map $\calP \to \calP_l$  yields a canonical
homomorphism $\pi \to \calP_l$.

\begin{proposition} \label{dense}
If $H_1(\pi,F)$ is finite dimensional, then each of the homomorphisms
$\pi \to
\calP_l$ has Zariski dense image.
\end{proposition}

\proof Denote the Zariski closure of the image of $\pi$ in $\calP_l$
by $\calZ_l$.  Each $\calZ_l$  is an algebraic subgroup of $\calP_l$.
Since
the composite
$$
H_1(\pi;F) \to H_1(\calZ_l) \to H_1(\calP_l)
$$
is an isomorphism, it follows that the second map is surjective.
Since
$\calP_l$ is unipotent, this implies that the inclusion
$\calZ_l \hookrightarrow \calP_l$ is surjective.  \endproof

\begin{theorem} \label{description}
If $H_1(\pi,F)$ is finite dimensional, then the natural map $\pi \to
\calP$ is the  Malcev completion of $\pi$ over $F$.
\end{theorem}

\proof  By the universal mapping property
(\ref{ump}), there is a canonical homomorphism from the
Malcev completion $\calU_F$ of $\pi$ to $\calP$.  It follows from
(\ref{dense}) that this homomorphism is surjective.

We now establish injectivity.  Suppose that $U$ is a unipotent
group defined over $F$.  This means that $U$ can be
represented as a subgroup of the group of upper triangular unipotent
matrices in $GL_n(F)$ for some $n$.  A representation $\rho : \pi \to
U$ induces a ring homomorphism
$$
\tilde{\rho} : F\pi \to gl_n(F).
$$
Since the representation is unipotent, it follows that
$\tilde{\rho}(J)$
is contained in the set of nilpotent upper triangular matrices. This
implies
that the
kernel of $\tilde{\rho}$ contains $J^n$. Consequently, $\rho$ induces
a
homomorphism
$$
\bar{\rho} : F\pi /J^n \to gl_n(F).
$$
Since the image of $J$ is contained in the nilpotent upper triangular
matrices, the image of the subgroup $\calP_{n-1}$ of $1+J/J^n$ is
contained in the group of unipotent upper triangular matrices.
Because
the image of $\pi$ is Zariski dense in $\calP_{n-1}$, it follows that
the image of $\calP_{n-1}$ is contained in $U$.  That is, there is a
homomorphism $\calP_{n-1} \to U$ of linear algebraic groups over
$F$ such that the diagram
$$
\matrix{
\pi &\to &\calP_{n-1} \cr
& \mapse{\rho} &\downarrow \cr
&& U \cr }
$$
commutes.  It follows that $\calU_F \to \calP$ is injective.
\endproof

\begin{corollary} \label{forms}  If $\calP_F$ is the Malcev
completion of the group $\pi$ over $F$, then the natural
homomorphism
$$
\calP_F \to \calP_\Q(F)
$$
is an isomorphism.
\end{corollary}

\proof This follows as $\calP_F$ is the set of group-like elements of
$F\pi\comp$, while  $\calP_{\Q}(F)$ is the set of group-like
elements of $\Q \pi \comp \otimes F\approx F\pi\comp$. \endproof

\section{Basic theory of relative completion}
\label{basic}

In this section we establish several basic properties of relative
completion. We begin by recalling
some basic facts about group extensions.  Once again, $F$ will denote
a  fixed field of characteristic zero.  All algebraic groups will be
linear.

Suppose that $L$ is an abstract group and that $A$ is an $L$ module.
The group of congruence classes of extensions
$$
0 \to A \to G \to L \to 1,
$$
where the natural action of $L$ on $A$ is the given action, is
naturally isomorphic to $H^2(L;A)$.  The identity is the semidirect
product $L \semi A$ (\cite[Theorem 4.1, p.~112]{maclane}). If
$H^1(L;A)$
vanishes, then any 2 splittings $s_0,s_1 : L \to L \semi A$ are
conjugate via an element of $A$ (\cite[Prop.\ 2.1,p.~106]{maclane}).
That is, there exists $a \in A$ such that $s_1 = a s_0 a^{-1}$.

If $S$ is a connected semisimple   algebraic group over $F$, and
if $A$ is a rational representation of $S$, then every extension
$$
0 \to A \to G \to S \to 1
$$
in the category of   algebraic groups over $F$ splits.  Moreover,
any 2 splittings $s_0,s_1 : S \to G$ are conjugate by an element of
$A$ \cite[p.~185]{humphreys}.

These results extend to the case when the kernel is unipotent.  For
this we need the following construction.

\begin{para} {\bf Construction.}
\label{construction}
Suppose that
\begin{equation}\label{G}
1 \to U \to G \to L \to 1
\end{equation}
is an extension of an abstract group $L$ by a group $U$.
Suppose that $Z$ is a central subgroup of $U$ and that the extension
\begin{equation}\label{G/Z}
1 \to U/Z \to G/Z \to L \to 1
\end{equation}
is split.  Let $s: L \to G/Z$ be a splitting.  Pulling back the
extension
$$
1 \to Z \to G \to G/Z \to 1
$$
along $s$, one obtains an extension
\begin{equation}\label{E_s}
1 \to Z \to E_s \to L \to 1.
\end{equation}
This determines, and is determined by, a class $\zeta_s$ in
$H^2(L;Z)$
which depends only on $s$ up to inner automorphisms by elements of
$G$.
The extension (\ref{E_s}) splits if and only if $\zeta_s = 0$.
\end{para}

\begin{proposition}\label{split}
If $\zeta_s = 0$, then the extension (\ref{G})
splits. Moreover, if any 2 splittings of (\ref{G/Z}) are conjugate,
and
if $H^1(L;Z) = 0$, then any 2 splittings of (\ref{G}) are conjugate.
\end{proposition}

\proof  If $\zeta_s$ vanishes, there is a splitting $\sigma: L \to
E_s$.
By composing $\sigma$ with the inclusion of $E_s \hookrightarrow G$,
one
obtains a splitting of the extension (\ref{G}).

Suppose now that $H^1(L;Z)$ vanishes, and that any 2 splittings
of $L \to G/Z$ are conjugate.  If $s_0,s_1 : L \to G$ are 2
splittings
of (\ref{G}), then their reductions $\s_0, \s_1 : L \to G/Z$ mod
$Z$ are conjugate. We may therefore assume that $s_0$ and $s_1$
agree mod $Z$.  The images of $s_0$ and $s_1$ are then contained
in the subgroup $E_s$ of $G$ determined by the section
$\s_0 = \s_1 : L \to G/Z$.   Since $H^1(L;Z) = 0$, there is an
element
$z$ of $Z$ which conjugates $s_0$ into $s_1$. \endproof
\medskip

A similar argument, combined with the facts about extensions of
algebraic groups at the beginning of this section, can be used
to prove the following result.

\begin{proposition}
Suppose that $S$ is a connected semisimple   algebraic group over
$F$.  If
$$
1 \to U \to G \to S \to 1
$$
is an extension in the category of   algebraic groups over $F$, and
if $U$ is
unipotent, then the extension splits, and any two
splittings $S \to G$ are conjugate. \endproof
\end{proposition}

By choosing compatible splittings and then taking inverse limits, we
obtain the analogous result for proalgebraic groups.

\begin{proposition}
Suppose that $S$ is a connected semisimple   algebraic group over
$F$.  If
$$
1 \to \calU \to \calG \to S \to 1
$$
is an extension in the category of   proalgebraic groups over $F$,
and
if $\calU$ is prounipotent, then the extension splits, and any two
splittings $S \to \calG$ are conjugate. \endproof
\end{proposition}

Suppose that $\Gamma$ is a abstract group, $S$ an
algebraic group over $F$, and that $\rho : \Gamma \to S$ a
representation with Zariski dense image.   Denote the image of
$\rho$ by $L$ and its kernel by $T$.  Thus we have an extension
$$
1 \to T \to \Gamma \to L \to 1.
$$
Let $\Gamma \to \calG_F$ be the completion over $F$ of $\Gamma$
relative to
$\rho$ and $\calU_F$ its prounipotent radical.   There is a
commutative diagram
$$
\matrix{
1 & \to & T & \to & \Gamma & \to & L & \to 1 \cr
&& \downarrow && \downarrow && \downarrow && \cr
1 & \to & \calU_F & \to & \calG_F & \to & S & \to 1 \cr}
$$
Denote the unipotent (i.e., the  Malcev) completion  of $T$ over $F$
by
$\calT_F$.  The universal mapping property
of $\calT_F$ gives a homomorphism $\Phi : \calT_F \to \calU_F$ of
prounipotent groups whose composition with the natural map $T\to
\calT_F$ is the  canonical map $T \to \calU_F$.  Denote the kernel of
$\Phi$ by $\calK_F$.

\begin{proposition}
Suppose that $H_1(T;F)$ is a finite dimensional.  If the action of
$L$
on $H_1(T;F)$ extends to a rational action of $S$, then the kernel
$\calK_F$ of $\Phi$ is central in $\calT_F$.
\end{proposition}

\proof First, $\Gamma$ acts on the
completion $FT\comp$ of the group algebra of $T$ by conjugation.
This action preserves the filtration by the powers of $\Jhat$, so it
acts on the associated graded algebra
$$
Gr^\bullet_J FT\comp = \bigoplus_{m=0}^\infty \Jhat^m/\Jhat^{m+1}.
$$
If $H_1(T;F)$ is finite dimensional, each truncation $FT/J^l$ of
$FT\comp$ is finite dimensional. This implies that each
of the groups $\Aut FT/J^l$ is an   algebraic group.  Since
$\Aut Gr^\bullet_J FT\comp$ is generated by $J/J^2 = H_1(T;F)$,
it follows that when $H_1(T;F)$ is finite dimensional, $\Aut
FT\comp$, the group of augmentation preserving algebra
automorphisms of $FT\comp$, is a proalgebraic group which is an
extension of a subgroup of $\Aut H_1(T;F)$ by a prounipotent group.

If the action of $\Gamma$ on $H_1(T;F)$ factors through a rational
representation $S \to \Aut H_1(T;F)$ of $S$, then we can form a
proalgebraic group extension
$$
1 \to J^{-1} \Aut FT\comp \to \calE \to S \to 1
$$
of $S$ by the prounipotent radical of $\Aut FT\comp$ by
pulling back the extension
$$
1 \to J^{-1} \Aut FT\comp \to \Aut FT\comp \to \Aut H_1(T;F)
$$
along $S\to \Aut H_1(T;F)$.  The representation
$\Gamma \to \Aut FT\comp$ lifts to a representation
$\Gamma \to \calE$ whose composition with the projection
$\calE\to S$ is $\rho : \Gamma \to S$.  This induces a
homomorphism $\calG_F \to \calE$.  Since the composite
$$
\calT_F \to \calG_F \to \calE \to \Aut FT\comp
$$
is the action of $\calT_F \subseteq FT\comp$ on $FT\comp$ by inner
automorphisms, it follows that the kernel of this map is the center
$Z(\calT_F)$ of $\calT_F$. It follows that the kernel of
$\calT_F \to \calU_F$ is a subgroup of $Z(\calT)$. \endproof
\medskip

The following result and its proof were communicated to me by
P.~Deligne \cite{deligne:letter}.

\begin{proposition}
Suppose that $S$ is semisimple and that  the natural action of $L$ on
$H_1(T;F)$ extends to a rational representation of $S$.  If $H^1(L;A)
= 0$ for
all rational representations $A$ of $S$, then $\Phi$ is surjective.
\end{proposition}

\proof The homomorphism $\calT_F \to \calU_F$ is surjective if and
only if the induced map $H_1(T;F) \to H_1(\calU_F)$ is surjective.
Let $A$ be the cokernel of this map.  This is a rational
representation
of $S$ as both $H_1(T;F)$ and $H_1(\calU_F)$ are.  By pushing out
the extension
$$
1 \to \calU_F\to \calG_F \to S \to 1
$$
along the map $\calU_F\to H_1(U) \to A$, we obtain an extension
$$
1 \to A \to G \to S \to 1
$$
of algebraic groups and a homomorphism $\Gamma \to G$ which lifts
$\rho$ and has
Zariski dense image in $G$.  Let
$$
1 \to A \to G^L \to L \to 1
$$
be the restriction of this extension to $L$.  Since $A$ is the
cokernel
of the map $\calT_F \to H_1(U)$, the image of $\Gamma$ in $G$ is
$\Gamma /T = L$. So $\rho$ induces a homomorphism $\rhotilde : L \to
G$ which has Zariski dense image.  The image of $\rhotilde$ lies in
$G^L$ and induces a splitting $\sigma : L \to G^L$ of the projection
$G^L \to L$. Since $G$ is an algebraic group, there is a splitting
$s : S \to G$ of the projection in the category of algebraic
groups.  This restricts to a splitting $s' : L \to G$.   Since
$H^1(L;A)$ vanishes, there exists $a\in A$ which conjugates
$s'$ into $\sigma$. Thus the image of $\sigma$ is contained in
the algebraic subgroup $as(S)a^{-1}$ of $G$.  Since the image of
$\sigma$ in $G$ is Zariski dense, it follows that $A$ must be
trivial.
\endproof
\medskip

A direct consequence of this result is a criterion for the map
$\rho : \Gamma \to S$ itself to be the $\rho$ completion. This
criterion is satisfied by arithmetic groups in semisimple groups
where each factor has $\Q$ rank $\ge 2$ \cite{ragunathan}.

\begin{corollary}
\label{comp_L}
Suppose that $\rho : \Gamma \to S$ is a homomorphism of an abstract
group into a semisimple algebraic group. If $H^1(L;A)$ vanishes for
all rational representations of $S$, and if $H^1(T;F)$ is zero (e.g.,
$\rho$ injective), then the relative completion of $\Gamma$ with
respect to $\rho$ is $\rho : \Gamma \to S$. \endproof
\end{corollary}

Next we consider the problem of imbedding an extension of $L$ by
a unipotent group $U$ in an extension of $S$ by $U$.

\begin{proposition}\label{lift}
Suppose that
$$
1 \to U \to G \to L \to 1
$$
is a split extension of abstract groups, where $U$ is a unipotent
group over $F$ and where the action of $L$ on $H^1(U)$ extends to
a rational representation of $S$.  If $H^1(L;A)$ vanishes for all
rational representations of $S$, then there is an extension
$$
1 \to U \to \tilde{G} \to S \to 1
$$
of algebraic groups
such that the original extension is the restriction of this one to
$L$.
\end{proposition}

\proof Since the first extension splits, we can write it as a semi
direct
product
$$
G = L \semi U.
$$
Denote the Lie algebra of $U$ by $\u$.  The group of Lie algebra
automorphisms of $\u$ is an algebraic group over $F$.  It can be
expressed
as an extension
$$
1 \to J^{-1} \Aut \u \to \Aut \u \to \Aut H_1(\u)
$$
where the kernel consists of those automorphisms which act trivially
on $H^1(\u)$, and therefore on the graded quotients of the lower
central
series of $\u$.  It is a unipotent group.  We can pull this extension
back along the representation $S \to \Aut H_1(\u)$ to
obtain an extension
$$
1 \to J^{-1} \Aut \u \to \tilde{A} \to S \to 1
$$
The representation $L \to \Aut \u$ lifts to a homomorphism
$L \to \tilde{A}$. This induces a homomorphism of the $S$-%
completion of $L$ into $\tilde A$.  By (\ref{comp_L}), the completion
of $L$ is simply the inclusion $L \hookrightarrow S$.  That is, the
representation $L \to
\tilde A$ extends to an algebraic group homomorphism $S \to \tilde
A$. This implies that
the representation of $L$ on $\u$ extends to a rational
representation
$S \to \Aut \u$.  We can therefore form the semi direct product
$S \semi U$, which is an algebraic group.  The homomorphism
$L \semi U \to S \semi U$ exists because of the compatibility of the
actions of $L$ and $S$ on $U$. \endproof
\medskip

Combining this with (\ref{split}), we obtain:

\begin{corollary}\label{split-lift}
Suppose that $\calU$ is a prounipotent group over $F$ with
$H_1(\calU)$
finite dimensional, and suppose that $\calZ$ is a central subgroup of
$\calU$. Suppose that
\begin{equation}\label{original}
1 \to \calU \to G \to L \to 1
\end{equation}
is an extension of abstract groups where the action of $L$ on
$H_1(\calU)$ extends to a rational representation of $S$.  Suppose
that
$$
1 \to \calU/\calZ \to \calG \to S \to 1
$$
is an extension  of proalgebraic groups which gives the extension
$$
1 \to \calU/\calZ \to G/\calZ \to L \to 1
$$
when restricted to $L$. If the class in $H^2(L;Z)$ given by
(\ref{construction}) vanishes, then there exists an extension
$$
1 \to \calU \to \tilde{\calG} \to S \to 1
$$
of proalgebraic groups whose restriction to $L$ is the extension
(\ref{original}). \endproof
\end{corollary}

By pushing out the extension
$$
1 \to T \to \Gamma \to L \to 1
$$
along the homomorphism $T \to \calT_F$, we obtain a ``fattening''
$\Gammahat$ of $\Gamma$.  Let $\calG^L$ be the inverse image of
$L$ in $\calG_F$. Using the universal mapping property of pushouts,
one can show easily that the natural homomorphism
$\Gamma \to \calG_F$ induces
homomorphism $\Gammahat \to \calG^L$.  These groups fit
into a commutative diagram of extensions:
$$
\matrix{
1 & \to & T & \to & \Gamma & \to & L & \to & 1 \cr
&& \downarrow && \downarrow && \| && \cr
1 & \to & \calT_F & \to & \Gammahat & \to & L & \to & 1  \cr
&& \downarrow && \downarrow && \| && \cr
1 & \to & \calU_F & \to & \calG^L & \to & L & \to & 1 \cr
&& \| && \downarrow && \downarrow && \cr
1 & \to & \calU_F & \to & \calG_F & \to & S & \to & 1 \cr}
$$

Next we introduce conditions we need to impose on our extension
for the remainder of the section.

\begin{para}  \label{conditions}
We will now assume that the extension $\Gamma$ of
$L$ by $T$ satisfies the following conditions. First, $H_1(T;F)$ is
finite dimensional and the action of $L$ on it extends to a rational
representation of $S$.   Next, we assume that $H^1(L;A)$ vanishes
for all rational representations of $S$. Finally, we add the new
condition that $H^2(L;A)$ vanishes for all {\it nontrivial}
irreducible
rational representations of $S$.  These conditions are satisfied by
arithmetic groups in semisimple groups where each factor has real
rank at least 8 \cite{borel:triv,borel:twisted}.
\end{para}

The following result is an immediate consequence of
(\ref{split-lift}).

\begin{proposition}\label{central}
If the conditions (\ref{conditions}) are satisfied, then
$\calK_F = \ker \Phi$ is contained in the center of the thickening
$\Gammahat$ of $\Gamma$.  \endproof
\end{proposition}

Applying the construction (\ref{construction}) to the thickening
$\Gammahat$ of $\Gamma$ and a splitting of $\calG^L \to L$, we
obtain an extension
\begin{equation}\label{obstn}
0 \to \calK_F \to G \to L \to 1
\end{equation}
which is unique up to isomorphism.  It follows from (\ref{central})
that this is a central extension.  Since $H_1(L;F)$ vanishes, there
is a
universal central extension with kernel an $F$
vector space.  It is the extension
$$
0 \to H_2(L;F) \to \tilde{L} \to L \to 1
$$
with cocycle the identity map
$$
\left\{H_2(L;F) \stackrel{id}{\to} H_2(L;F)\right\}\in
\Hom(H_2(L;F),H_2(L;F)) \approx H^2(L;H_2(L;F)).
$$
The central extension (\ref{obstn}) is classified by a linear map
$\psi_F : H_2(L;F) \to \calK_F$.  Because all splittings
$s: L \to \calG^L$ are conjugate (\ref{split}), the class of this
extension is independent of the choice of the splitting.

Since $\calK_F$ is an abelian unipotent group,
$\calK_F(k) = \calK_F\otimes k$ for all fields $k$ which contain $F$.
The homomorphism $\psi_F$ satisfies the
following naturality property:

\begin{proposition}
If $k$ is an extension field of $F$, then the diagram
$$
\matrix{
H_2(L;k) &\to & \calK_k \cr
|| && \downarrow \cr
H_2(L;F)\otimes k & \to & \calK_F(k) \cr }
$$
commutes. \endproof
\end{proposition}

The next result bounds the size of $\calK_F$.

\begin{proposition}\label{surj}
If the conditions (\ref{conditions}) hold, then the natural map
$\psi_F : H_2(L;F) \to \calK_F$ is surjective.
\end{proposition}

\proof As above, we shall denote by $G$ the central extension of $L$
by
$\calK_F$.  Let $A$ be the cokernel of $\psi_F$ and $E$ the cokernel
of $\psi_F : H_2(L;F) \to G$.  Then $E$ is a central extension of $L$
by $A$. Because the composite $H_2(L;F) \to \calK_F \to A$ is
trivial,
this extension is split.  From (\ref{split}) it follows
that the extension
$$
1 \to \calT_F/\im \psi_F \to \Gammahat/\im \psi_F \to L \to 1
$$
is split. By (\ref{lift}), this implies that there is a proalgebraic
group
$\calE$ which is a semidirect product of $S$ by $\calT_F/\im \psi_F$
into which $\Gammahat/\im\psi_F$ injects.  The map of
$\Gamma$ to $\calE$ induces a map $\calG_F \to \calE$.  Since the
kernel
of the map $\calT_F \to \calE$ is $\im \psi_F$, it follows that
$\calK_F$
is contained in $\im \psi_F$. \endproof

\begin{corollary}\label{indep}
If the conditions (\ref{conditions}) hold, then the natural map
$\calG_k \to \calG_F(k)$ associated to a field extension $k:F$ is an
isomorphism.
\end{corollary}

\proof   By (\ref{forms}), the natural
map $\calT_F \to \calT_F(k)$ is an isomorphism.  Since $\calT_K \to
\calU_K$ is surjective for all fields $K$, the natural map $\calU_k
\to
\calU_F(k)$ is also surjective.  Consequently, the natural map $K_k
\to K_F(k)$ is injective. Since $\calK_F$ is abelian unipotent,
$\calK_F(k) = \calK_F \otimes k$. But it follows from (\ref{surj})
that
$\calK_k \to \calK_F(k)$ is surjective, and therefore an isomorphism.
\endproof

\section{The Johnson homomorphism}
\label{johnson_homo}

This section is a brief review of the construction of Johnson's
homomorphism \cite{johnson:survey,johnson:2}.  There are two
equivalent ways,  both due to Johnson, to define a homomorphism
$$
T_g \to \Lambda^3 H_1(C;\Z)
$$
where $C$ is a compact Riemann surface of genus $g$.

Choose a base point $x$ of $C$.  The first construction uses the
action
of $T_g^1$ on $\pi_1(C,x)$. Denote the lower central series
of $\pi_1(C,x)$ by
$$
\pi_1(C,x) = \pi^1 \supseteq \pi^2 \supseteq \pi^3 \supset \cdots
$$
The first graded quotient $\pi^1/\pi^2$ is $H_1(C;\Z)$. The second is
naturally isomorphic to $\Lambda^2 H_1(C;\Z)/\langle q\rangle$,
where
$q : \Lambda^2 H^1(C;\Z) \to H^2(C;\Z) \approx\Z$ is the cup product.
If $\gamma_j$, $j=1,\ldots,2g$, are generators of $\pi_1(C,x)$, then
the residue class of the commutator
$\gamma_j\gamma_k\gamma_j^{-1}\gamma_k^{-1}$ modulo
$\pi^3$ is the element $c_j\wedge c_k$ of
$\Lambda^2 H_1(C;\Z)/\langle q \rangle$,
where $c_j$ denotes the homology class of $\gamma_j$.  The  form
$q$
is just the equivalence class of the standard relation in
$\pi_1(C,x)$.

If $\phi : (C,x) \to (C,x)$ is a diffeomorphism which represents an
element of $T_g^1$, then $\phi$ acts trivially on the homology of
$C$.
It therefore acts as the identity on each graded quotient of the
lower
central series of $\pi_1(C,x)$.  Define a function  $\pi \to \pi$ by
taking $\gamma$ to $\phi(\gamma)\gamma^{-1}$.  Since $\phi$ acts
trivially on $H_1(C)$, it follows that this map takes $\pi^l$ into
$\pi^{l+1}$.  In particular, it induces a well defined function
$$
\tilde{\tau}(\phi) : H_1(C;\Z) \to
\Lambda^2 H_1(C;\Z)/\langle q \rangle
$$
between the first two graded quotients of $\pi$, which is easily seen
to
be linear.  Using Poincar\'e
duality, $\tilde{\tau}(\phi)$ can be regarded as an element of
$$
H_1(C;\Z) \otimes
\left(\Lambda^2 H_1(C;\Z)/\langle q \rangle\right).
$$
The map $\phi \mapsto \tilde{\tau}(\phi)$ induces a group
homomorphism
$$
T_g^1 \to H_1(C;\Z) \otimes
\left(\Lambda^2 H_1(C;\Z)/\langle q \rangle\right).
$$
and therefore a homomorphism
$$
\hat{\tau} : H_1(T_g^1) \to H_1(C;\Z) \otimes
\left(\Lambda^2 H_1(C;\Z)/\langle q \rangle\right).
$$
There is a natural inclusion
$$
\Lambda^3 H_1(C;\Z) \to  H_1(C;\Z) \otimes
\left(\Lambda^2 H_1(C;\Z)/\langle q \rangle\right).
$$
defined by
$$
x\wedge y\wedge z \mapsto x\otimes (y\wedge z) +
y \otimes (z \wedge x) + z \otimes (x \wedge y).
$$
Johnson has proved that the  image of $\hat{\tau}$ is contained in
the image of this map, so that $\hat{\tau}$ induces a homomorphism
$$
\tau_g^1 : H_1(T_g^1) \to \Lambda^3 H^1(C;\Z).
$$
It is not difficult to check that this homomorphism is $Sp_g(\Z)$
equivariant.

This story can be extended to $T_g$ as follows.  There is a natural
extension
$$
1 \to \pi_1(C,x) \to T_g^1 \to T_g \to 1.
$$
Applying $H_1$, we obtain the diagram
$$
\matrix{
H_1(C;\Z) & \to & H_1(T_g^1) &\to & H_1(T_g) &\to & 0\cr
&&\mapdown{\tau_g^1} && &&\cr
&&\Lambda^3 H_1(C;\Z)&&&&\cr}
$$
whose top row is exact.
Identify the lower group with $H_3(\Jac C;\Z)$.
Since $\Jac C$ is a group with torsion free homology, the group
multiplication induces a product on its homology which is
called the {\it Pontrjagin product}. The composite of
$\tau_g^1$ with the map from $H_1(C)$ takes a class in $H_1(\Jac C)$
to its Pontrjagin
product with $q =[C] \in H_2(\Jac C;\Z)$.  It follows that $\tau_g^1$
induces a map
$$
\tau_g : H_1(T_g) \to
\Lambda^3 H_1(C;\Z)/\left(q\wedge H_1(C;\Z)\right).
$$

The following fundamental theorem is due to Dennis Johnson.

\begin{theorem}\label{johnson:3}
When $g\ge 3$,
the homomorphisms $\tau_g^1$ and $\tau_g$ are isomorphisms
modulo 2 torsion.
\end{theorem}

The second way to construct a map  $T_g^1 \to \Lambda^3 H_1(C;\Z)$
is as follows.  Suppose that $\phi$ represents an element of $T_g^1$.
Let $M_\phi \to S^1$ be the bundle over the circle constructed by
identifying the point $(z,1)$ of $C\times [0,1]$ with the point
$(\phi(z),0)$.  Since $\phi$ fixes the basepoint $x$, the map
$t \to (x,t)$ induces a section of $M_\phi \to S^1$.  Similarly, one
can
construct the bundle of jacobians; this is trivial as $\phi$ acts
trivially  on $H_1(C)$.  One can imbed $M_\phi$ in this bundle of
jacobians using this section of basepoints. Let $p$ be the
projection of the bundle of jacobians onto one of its fibers.  Then
one has the 3 cycle $p_\ast M_\phi$ in $\Jac C$.

\begin{proposition}\label{tau_2}
{\rm \bf \cite{johnson:survey}}
When $g\ge 3$, the homology class
$$
p_\ast [M_\phi] \in H_3(\Jac C;\Z) \approx \Lambda^3 H_1(C;\Z)
$$
is $\tau_g^1(\phi)$. \endproof
\end{proposition}

This can be proved, for example, by checking that both maps agree
on what Johnson calls ``bounding pair" maps. These maps generate
the Torelli group when $g\ge 3$ \cite{johnson:1}.

\section{The cycle $C - C^-$}
\label{cycle}

In this section we relate the algebraic cycle $C - C^-$ to Johnson's
homomorphism.

Let $C$ be a compact Riemann surface. Denote its jacobian by
$\Jac C$.  This is defined to be $\Pic^0 C$, the group of divisors of
degree zero on $C$ modulo principal divisors.   Each divisor $D$ of
degree zero may be written as the boundary of a topological 1-chain:
$D = \partial \gamma$.  Taking $D$ to the functional
$\omega \mapsto \int_\gamma \omega$ on the space $\Omega(C)$
of holomorphic 1-forms yields a well defined map
\begin{equation}\label{AJ}
\Pic^0 C \to \Omega(C)^\ast/H_1(C;\Z).
\end{equation}
This is an isomorphism by Abel's Theorem.

For each $x\in C$, we have an Abel-Jacobi map
$$
\nu_x : C \to \Jac C
$$
which is defined by $\nu_x(y) = y - x$.   This map is an imbedding if
the genus $g$ of $C$ is $\ge 1$.  Denote its image by $C_x$. This is
an
algebraic 1-cycle in $\Jac C$.  Denote the involution $D\mapsto -D$
of
$\Jac C$ by $i$, and the image of $C_x$ under this involution
by $C_x^-$.  Since $i$ induces $-id$ on $H^1(\Jac C;\Z)$, and since
$i^\ast$ is a ring homomorphism, we see that $i^\ast$ acts as
$(-1)^k$
on $H_k(\Jac C)$.  It follows that $C_x$ and $C_x^-$ are
homologically
equivalent.  Set $Z_x = C_x - C_x^-$.  This is a homologically
trivial
1-cycle.

Griffiths has a construction which associates to a homologically
trivial
analytic cycle in a compact K\"ahler manifold a point in a complex
torus. His construction is a generalization the construction of
the map (\ref{AJ}).   We review this construction briefly.
Suppose that $X$ is a compact K\"ahler manifold, and that $Z$ is an
analytic $k$-cycle in $X$ which is homologous to 0.   Write $Z =
\partial
\Gamma$, where $\Gamma$ is a topological $2k+1$ chain in $X$.
Define
$$
F^p H^m(X) = \bigoplus_{s \ge p} H^{s,m-s}(X).
$$
Each class in $F^pH^m(X)$ can be represented by a closed form where
each term of a local expression in terms of local holomorphic
coordinates $(z_1,\ldots, z_n)$ has at least $p$ $dz_j$s.
Integrating such representatives of classes in $F^{k+1}H^{2k+1}(X)$
over $\Gamma$ gives a well defined functional
$$
\int_\Gamma : F^{k+1} H^{2k+1}(X) \to \C.
$$
The choice of $\Gamma$ is unique up to a topological $2k+1$ cycle.
So $Z$ determines a point of the complex torus
$$
J_k(X) := F^{k+1} H^{2k+1}(X)^\ast/H_{2k+1}(X;\Z).
$$
This group is called the {\it $k$th intermediate jacobian of $X$}.

In our case, the cycle $Z_x$ determines a point $\zeta_x(C)$ in the
intermediate jacobian
$$
J_1(\Jac C) = F^2H^3(\Jac C)^\ast/H_3(\Jac C; \Z).
$$

The
homology class of $C_x$ is easily seen to be independent of $x$.
Taking
the Pontrjagin product with $[C]$ defines an injective map
$$
H_1(\Jac C;\Z) \hookrightarrow H_3(\Jac C;\Z).
$$
The dual $H^3(\Jac C) \to H^1(C)$ is a morphism of Hodge structures
of type $(-1,-1)$ --- that is, $H^{s,t}$ gets mapped into
$H^{s-1,t-1}$.
It follows that this map induces an imbedding
$$
\Phi : \Jac C \hookrightarrow J_1(\Jac(C))
$$
of complex tori.
We shall denote the cokernel of this map by $JQ_1(\Jac C)$. It is
trivial when $g < 3$.

The following result is not difficult; a proof may be found in
\cite{pulte}.

\begin{proposition}
If $x,y \in C$, then
$$
\zeta_x(C) - \zeta_y(C) = 2 \Phi(x-y).
$$
In particular, the image of $\zeta_x(C)$ in $JQ_1(\Jac C)$ is
independent of $x$.
\end{proposition}

Denote the common image of the $\zeta_x(C)$ in
$JQ_1(\Jac C)$ by $\zeta(C)$.

We now suppose that the genus $g$ of $C$ is $\ge 3$.
Fix a level $l$ so that the moduli space $\M_g^n(l)$ of curves of
genus $g$ and $n$ marked points with a level $l$ structure is
smooth. (Any
$l\ge 3$ will do for the time being.)  Denote the space of
principally
polarized abelian varieties of dimension $g$ with a level $l$
structure
by $\A_g(l)$.

One can easily construct  bundles $\calJ_m \to \A_g(l)$ of complex
tori
whose fiber over the abelian variety $A$ is $J_m(A)$. The pullback of
$\calJ_0$ along the period map $\M_g(l) \to \A_g(l)$
is the bundle of jacobians associated to the universal curve.  The
imbedding $\calJ_0 \hookrightarrow \calJ_1$ defined over $\M_g(l)$
extends to all of $\A_g(l)$ as the homology class $[C] \in H_2(\Jac
C;\Z)$
extends to a class $q \in H_2(A;\Z)$ for every abelian variety
$A$. Denote the bundle of quotient tori by $\calQ \to \A_g(l)$.

Denote the level $l$ congruence subgroup of $Sp_g(\Z)$ by $L(l)$.
Since $\A_g(l)$ is an Eilenberg-Mac Lane space $K(L(l),1)$, and since
the fiber of $\calJ_1 \to \A_g(l)$ is a $K(H_3(\Jac C ;\Z),1)$, it
follows
that $\calJ_1$ is also an Eilenberg-Mac Lane space whose
fundamental group is an extension of $L(l)$ by $H_3(\Jac C;\Z)$.
Since this bundle has a section, viz., the zero section, this
extension splits.  The action of $L(l)$ on $H_3(\Jac C;\Z)$ is the
restriction of  the third exterior power of the fundamental
representation of $Sp_g$.  There is a similar story with $\calJ_1$
replaced by $\calQ$.

We shall denote the quotient
$\Lambda^3 H_1(C;\Z)/\left([C]\cdot H_1(C;\Z)\right)$ by
$Q\Lambda^3 H_1(C)$.  This is the fundamental group
of $JQ_1(\Jac C)$.

\begin{proposition}
The spaces $\calJ_1$ and $\calQ$ are Eilenberg-Mac Lane spaces with
fundamental groups
$$
\pi_1(\calJ_1,0_C) \approx L(l) \semi \Lambda^3 H_1(C;\Z)
$$
and
$$
\pi_1(\calQ,0_C) \approx
L(l) \semi Q\Lambda^3 H_1(C),
$$
respectively.
Here $0_C$ denotes the identity element in the fiber over $\Jac C$.
\end{proposition}

The normal functions $(C,x) \mapsto \zeta_x(C)$ and $C \mapsto
\zeta(C)$ give
lifts of the period map:
$$
\matrix{
&& \calJ_1 \cr
& \mapne{\zeta_g^1} & \downarrow \cr
\M_g^1(l) & \to & \A_g(l) \cr }
\qquad
\matrix{
&& \calQ \cr
& \mapne{\zeta_g} & \downarrow \cr
\M_g(l) & \to & \A_g(l) \cr }
$$
These induce maps of fundamental groups
$$
{\zeta_g^1}_\ast : \Gamma_g^1(l) \to L(l)\semi \Lambda_3 H_1(C)
$$
and
$$
{\zeta_g}_\ast : \Gamma_g(l) \to L(l) \semi Q\Lambda^3 H_1(C)
$$
Since these maps commute with the canonical projections to $L(l)$,
these induce $L(l)$ equivariant maps
$$
\zeta_g^1 : H_1(T_g^1) \to \Lambda^3 H_1(C;Z), \quad
\zeta_g^1 : H_1(T_g) \to Q\Lambda^3 H_1(C;Z)
$$

The following result follows easily from (\ref{tau_2}). The factor
of 2 arises as both $C$ and $C^-$ each contribute a copy of the
Johnson homomorphism.

\begin{proposition}\label{double}
The map $\zeta_g^n$ is twice Johnson's map $\tau_g^n$ for $n=0,1$.
\endproof
\end{proposition}

It is natural to try to give a ``motivic description'' of the Johnson
homomorphism rather than of twice it. Looijenga (unpublished) has
done this by constructing a normal function which compares the
cycle $C_x$ to a fixed topological (but not algebraic) cycle in $\Jac
C$ which
is homologous to $C_x$.  At the cost
of being more abstract,  we give another description
which does not make use of Looijenga's topological cycle.

We will only consider the pointed case, the unpointed case being
similar.  For each pointed curve $(C,x)$, the cycle $C_x$
determines a point $c_x$ in the Deligne cohomology group
$H^{2g-2}_\calD(\Jac(C),\Z(g-1))$.  This group is an extension
of the Hodge classes $H^{g-1,g-1}_\Z(\Jac C)$ by the intermediate
jacobian $J_1(\Jac C)$:
$$
0 \to J_1(\Jac C) \to H^{2g-2}_\calD(\Jac(C),\Z(g-1)) \to
H^{g-1,g-1}_\Z(\Jac C) \to 0
$$
One can consider the bundle over $\A_g(l)$ whose fiber over the
abelian variety $A$ is $H^{2g-2}_\calD(A,\Z(g-1))$.  The subbundle
whose fiber over $A$ is the $J_1(A)$ coset of the class of the
polarization $q \in H^{g-1,g-1}_\Z(A)$ is a principal $\calJ_1$
torsor
over $\A_g(l)$. Denote it by $\calZ \to \A_g(l)$.  The cycle gives
a lift
$$
\matrix{
&& \calZ \cr
& \mapne{c} & \downarrow \cr
\M_g^1(l) & \to & \A_g(l) \cr}
$$
of the period map. The total space $\calZ$ is an Eilenberg Mac Lane
space whose fundamental group is an extension of $L(l)$ by
$\Lambda^3 H_1(C;Z)$. The map induced by $c$ on fundamental groups
induces the Johnson homomorphism $\tau_g^1$. This follows
directly from (\ref{tau_2}).

\section{Completion of mapping class groups}
%\label{rel_comp_Mg}

Denote the completion of the mapping class group $\Gamma_{g,r}^n$
with respect to the canonical representation $\Gamma_{g,r}^n \to
Sp_g$ by $\calG_{g,r}^n$ and its prounipotent radical
by $\calU_{g,r}^n$. Denote the Malcev completion
of the Torelli group $T_{g,r}^n$ by $\calT_{g,r}^n$ and the kernel
of the natural homomorphism $\calT_{g,r}^n \to \calU_{g,r}^n$ by
$\calK_{g,r}^n$.  These groups are all defined over $\Q$
by (\ref{indep}).

For all $g \ge 3$, and all arithmetic subgroups $L$ of $Sp_g(\Z)$,
$H^1(L;A)$ vanishes for all rational representations $A$ of $Sp_g$.
By the results of Borel \cite{borel:triv,borel:twisted}, the
hypotheses
(\ref{conditions}) are satisfied by all arithmetic subgroups $L$ of
$Sp_g(\Z)$ and $H^2(L;\Q)$ is one dimensional when $g\ge 8$.
This, combined with (\ref{surj}) yields the following result.

\begin{proposition} For all $g \ge 3$ and all $n, r \ge 0$, the
natural
map $\calT_{g,r}^n \to \calU_{g,r}^n$ is surjective.  When $g\ge 8$,
the kernel $\calK_{g,r}^n$ is either trivial or isomorphic to $\Q$.
\endproof
\end{proposition}

Our main result is:

\begin{theorem}\label{main}
For all $g \ge 3$ and all $r,n\ge 0$, the kernel $\calK_{g,r}^n$ is
non-trivial, so that $\calK_{g,r}^n \approx \Q$ when $g\ge 8$.
\end{theorem}

Let $\lambda_1,\ldots,\lambda_g$ be a fundamental set of weights
of $Sp_g$.  For a dominant integral weight $\lambda$, denote the
irreducible representation with highest weight $\lambda$ by
$V(\lambda)$.

In this section we will reduce the proof of Theorem \ref{main} to the
proof of the the case $g=3$ and $r=n=0$.   The following assertion
follows easily by induction on $n$ and $r$ from Johnson's result.

\begin{proposition}
If $g \ge 3$, then for all $n,r\ge 0$, there is an $Sp_g$ equivariant
isomorphism
$$
H_1(T_{g,r}^n;\Q) \approx V(\lambda_3) \oplus V(\lambda_1)^{n+r}.
\quad\square
$$
\end{proposition}

The important fact for us is that the multiplicity of $\lambda_3$ in
$H_1(T_{g,r}^n;\Q)$ is
always 1.  Denote the second graded quotient of the lower central
series of $\calT_{g,r}^n$ by $\calV_{g,r}^n$.  The commutator induces
a linear surjection
$$
\Lambda^2 H_1(T_{g,r}^n;\Q) \to \calV_{g,r}^n
$$
which is $Sp_g$ equivariant.  By Schur's lemma, there is a unique
copy of the trivial representation in $\Lambda^2 V(\lambda_3)$.
Let $\beta_{g,r}^n : \Q \to \calV$ be the composite
$$
\Q \hookrightarrow \Lambda^2 V(\lambda_3) \to
\Lambda^2 H_1(T_{g,r}^n;\Q) \to \calV_{g,r}^n,
$$
where the first map is the inclusion of the trivial representation.
Consider the following assertions:
\medskip

\noindent{$\boldmath A_{g,r}^n$:} \qquad The map $\beta_{g,r}^n$ is
injective.
\medskip

\begin{proposition}\label{red_1}
If $h\ge g\ge 3$, $s\ge r\ge 0$ and $m\ge n\ge $, then $A_{g,r}^n$
implies $A_{h,s}^m$. Furthermore $A_{g,1}$ implies $A_g$.  In
particular, $A_3$ implies $A_{g,r}^n$ for all
$g \ge 3$ and $n,r \ge 0$.
\end{proposition}

\proof It is easy to see that $\beta_{g,r}^n$ is the composite of
$\beta_{g,s}^m$ with the canonical quotient map
$\calV_{g,s}^m \to \calV_{g,r}^n$ whenever $m\ge n$ and $s\ge r$. So
$A_{g,r}^n$ implies $A_{g,s}^m$.
Moreover, when $r\ge 1$, the composition of $\beta_{g,r}^n$ with the
map $\calV_{g,r}^n \to \calV_{g+1,r}^n$ induced by any one of the
natural maps $\Gamma_{g,r}^n \to \Gamma_{g+1,r}^n$ is
$\beta_{g+1,r}^n$.  So, in this case, $A_{g,r}^n$ implies
$A_{g+1,r}^n$.

To see that $A_{g,1}$ implies $A_g$,  consider the
group extension
$$
1 \to \pi_1(T^\ast_1 C,\vec{v}) \to T_{g,1} \to T_g \to 1
$$
where $T^\ast_1 C$ denotes the unit cotangent bundle of $C$. It is
not difficult to show that  $\calV_{g,1}$ is the direct sum of
$\calV_g$ and the second graded quotient of the lower central series
of $\pi_1(T^\ast_1 C,\vec{v})$, from which the assertion follows.
\endproof
\medskip

The assertions $A_{g,r}^n$ can be proved using
Harer's computations of $H^2(\Gamma_{g,r}^n;\Q)$
\cite{harer:h2,harer:h3}.  However, we will prove $A_3$ in
the course establishing the rest of Theorem \ref{main} directly,
without appeal to Harer's computation.

Denote the fattening (see discussion following (\ref{split})) of the
mapping class group $\Gamma_{g,r}^n$
by $\Gammahat_{g,r}^n$. This is an extension
$$
1 \to \calT_{g,r}^n \to \Gammahat_{g,r}^n \to Sp_g(\Z) \to 1.
$$
Dividing out by the commutator subgroup of $\calT_{g,r}^n$, we obtain
an
extension
\begin{equation}\label{aux}
0 \to H_1(T_{g,r}^n;\Q) \to \calE_{g,r}^n \to Sp_g(\Z) \to 1.
\end{equation}
This extension is split. This can be seen using a straightforward
generalization of (\ref{double}).

Now suppose $A_{g,r}^n$ holds.  Then there is a
quotient $G_{g,r}^n$ of $\calT_{g,r}^n$ which is an extension of
$H_1(T_{g,r}^n;\Q)$ by $\Q$. The cocycle of the extension being a
non-zero multiple of the polarization
$$
\theta \in \Lambda^2 V(\lambda_3) \subseteq H^2(T_{g,r}^n;\Q).
$$
Dividing $\Gammahat_{g,r}^n$ by the kernel of
$\calT_{g,r}^n \to G_{g,r}^n$, we obtain an extension
\begin{equation}\label{ext}
1 \to G_{g,r}^n \to E_{g,r}^n \to Sp_g(\Z) \to 1.
\end{equation}
Since the extension (\ref{aux}) splits, we can apply the construction
(\ref{construction}) to obtain an extension
$$
0 \to \Q \to H_{g,r}^n \to Sp_g(\Z) \to 1
$$
or equivalently, a class $e_{g,r}^n \in H^2(Sp_g(\Z);\Q)$.
By (\ref{split}), the non-triviality of this class is the obstruction
to
splitting the extension (\ref{ext}), which, by (\ref{lift}) is the
obstruction to imbedding it in an algebraic group extension of $Sp_g$
by $G_{g,r}^n$.  This proves the following statement.

\begin{proposition}
If $A_{g,r}^n$ holds and if $e_{g,r}^n$ is non-zero, then
$\calK_{g,r}^n$ is non-trivial.
\end{proposition}

To reduce the proof of Theorem (\ref{main}) to the genus 3 case, we
need to relate the classes $e_{g,r}^n$.

\begin{proposition}\label{red_2}
For fixed $g \ge 3$, the classes $e_{g,r}^n$ are all equal. The image
of $e_{g+1,1}$ under the natural map
$H^2(Sp_{g+1}(\Z);\Q) \to H^2(Sp_g(\Z);\Q)$ is $e_{g,1}$.
\end{proposition}

\proof Both statements follow from the naturality of the
construction. \endproof
\medskip

Combining (\ref{red_1}) and (\ref{red_2}), we have:

\begin{proposition}\label{main_3}
If $A_3$ holds and $e_3$ is non trivial, then Theorem \ref{main}
is true.
\end{proposition}

\section{Proof of Theorem \protect\ref{main}}
\label{proof}

We prove Theorem \ref{main} by proving (\ref{main_3}).
In this section we assume that the reader is familiar with mixed
Hodge theory.  We will use the notation and conventions of
\cite[\S\S 2--3]{hain:heights}.  The moduli space $\A_g$ can
be thought of
as the moduli space of principally polarized Hodge structures of
weight $-1$, level 1 and dimension $2g$; the abelian variety
$A \in \A_g$ corresponds to
the Hodge structure $H_1(A)$ and its natural polarization.
We can construct bundles over $\A_g$ by considering moduli
spaces of various mixed Hodge structures derived from such a Hodge
structure of weight $-1$.  To guarantee that we have a smooth
moduli space, we fix a level $l$ so that $\A_g(l)$ is smooth.

For a Hodge structure $H\in\A_g$, where $g \ge 3$,  with principal
polarization $q$, we define $QH$ to be the Hodge structure
$$
\left[\Lambda^3 H/\left( q\wedge H\right)\right] \otimes \Z(-1)
$$
which is of weight $-1$.  Denote the dual Hodge structure
$$
\Hom(QH,\Z(1)) \approx \
\ker\left\{ \_\wedge q : \Lambda^{2g-3} H \to
\Lambda ^{2g-1} H\right\}\otimes \Z(2-g)
$$
by $PH$.

The set of all mixed Hodge structures with weight graded quotients
$\Z$ and $QH$ is naturally isomorphic to the complex torus
$$
J(QH) := QH_\C / \left(F^0 QH_\C + QH_\Z \right).
$$
If $H = H_1(\Jac C)$, then $J(QH)$ is the torus $JQ_1(\Jac C)$
defined
in \S\ref{cycle}.
As we let $H$ vary over $\A_g(l)$, we obtain the bundle
$\calQ \to \A_g(l)$ of intermediate jacobians constructed in
\S\ref{cycle}.

The set of all mixed Hodge structures with weight graded quotients
$QH$ and $\Z(1)$ is naturally isomorphic to the complex torus
$$
J(PH) := PH_\C / \left(F^0 PH_\C + PH_\Z \right).
$$
This torus is the dual of $J(QH)$.
Performing this construction for each $H$ in $\A_g(l)$, we
obtain a bundle  $\calP \to \A_g(l)$ of complex tori.

We now construct a line bundle over the fibered product
$$
\calP\times_{\A_g(l)} \calQ \to \A_g(l).
$$
It is the {\it biextension line bundle}.  For a Hodge structure $H\in
\A_g$,
let $B(H)$ be the set of mixed Hodge structure with
weight graded quotients canonically isomorphic to $\Z$, $QH$, and
$\Z(1)$.  Set
$$
G_\Z = \pmatrix{
1 & QH_\Z & \Z(1) \cr
0 & 1 & PH_\Z \cr
0 & 0 & 1 \cr }
\quad
G = \pmatrix{
1 & QH_\C & \C \cr
0 & 1 & PH_\C \cr
0 & 0 & 1 \cr }
$$ $$
F^0 G_\Z = \pmatrix{
1 & F^0 QH & 0 \cr
0 & 1 & F^0 PH \cr
0 & 0 & 1 \cr }.
$$
There is a natural isomorphism
$$
B(H) = G_\Z\backslash G / F^0 G.
$$
The natural projection
$$
B(H) \to J(QH) \times J(PH),
$$
which takes $V \in B(H)$ to
$(V/W_{-2},W_{-1}V)$, is a principal $\C^\ast$ bundle.  Doing this
construction over $\A_g(l)$, we obtain a $\C^\ast$ bundle
$$
\calB \to \calQ\times_{\A_g(l)} \calP
$$
Denote the corresponding line bundle by $\calL \to
\calQ\times_{\A_g(l)}
\calP$.

Since the fiber of $\calB \to \A_g(l)$ over $H$ is the nilmanifold
$B(H)$, which is an Eilenberg-Mac Lane space, it follows that $\calB$
is
also an Eilenberg-Mac Lane space whose
fundamental group of $\calB$ is an extension
\begin{equation}\label{extn}
1 \to G_\Z \to \pi_1(\calB,\ast) \to L(l) \to 1.
\end{equation}
Taking $H\in \A_g(l)$ to the split biextension $\Z \oplus H \oplus
\Z(1)$ defines a section of the bundle $\calB \to \A_g$.  It follows
that the extension (\ref{extn}) is split.

\begin{proposition}
$\pi_1(\calB,\ast) \approx L(l) \semi G_\Z\quad $. \endproof
\end{proposition}

We now restrict ourselves to genus 3
and consider the problem of lifting the period map
$\M_3(l) \to \A_3(l)$ to $\calB$.  The idea is that such a lifting
will
be the period map of a variation of mixed Hodge structure.  To this
end we define certain algebraic cycles which are just more canonical
versions of the cycle $C-C^-$ considered in \S\ref{cycle}.

Suppose that $C$ is a curve of genus 3 and that $\alpha$ is a
theta characteristic of $C$ --- i.e., $\alpha$ is a square root of
the
canonical divisor $\kappa_C$.  Denote the algebraic cycle
corresponding
to the canonical inclusion $C \hookrightarrow \Pic^1C$ by $C$.
Denote the involution $x \mapsto \alpha -x$ of $\Pic^1 C$ by
$i^\alpha$.   For $D \in \Pic^0C$, denote the translation map
$x \mapsto x + D$ by
$$
\tau^D : \Pic^1 C \to \Pic^1 C.
$$
For $D \in \Pic^0 C$, define $C_D = \tau^D_\ast C$ and
$Z_{\alpha,D} = C_D - i^\alpha_\ast C_D$.
Each $Z_{\alpha,D}$ is homologous to zero.
Set $Z_\alpha = Z_{\alpha,0}$.

Let $\Theta_\alpha$ be the {\it theta divisor}
$$
\left\{ x+y - \alpha  : x,y \in C \right\} \subseteq \Pic^0 C
$$
and $\Delta$ be the {\it difference divisor}
$$
\left\{x - y : x, y \in C \right\} \subseteq \Pic^0 C.
$$
The following fact is easily verified.

\begin{proposition}\label{disjoint}
The cycles $Z_\alpha$ and $Z_{\alpha,D}$ have disjoint supports
if and only if $D \not\in \Theta_\alpha \cup \Delta$. \endproof
\end{proposition}

Now choose a point $\delta \in \Pic^0 C$ of order 2 such that
$\beta := \alpha + \delta$ is an {\it even} theta characteristic.
(i.e.,
$h^0(C,\beta) = 0$ or 2.)

\begin{proposition}\label{hyperelliptic}
The cycles $Z_\alpha$ and $Z_{\alpha,\delta}$ have disjoint supports
except when $C$ is hyperelliptic and either $\delta$ is the
difference
between two distinct Weierstrass points or $\alpha + \delta$ is the
hyperelliptic series.
\end{proposition}

\proof  By (\ref{disjoint}), $Z_\alpha$ and $Z_{\alpha,\delta}$
intersect if and only if $\delta \in \Delta$ or
$\delta \in \vartheta_\alpha$.  In the first  case, there exist
$x,y \in C$ such that $x-y = \delta\neq 0$. So $2x - 2y=0$, which
implies that $C$ is hyperelliptic and that $x$ and $y$ are distinct
Weierstrass points.   In the second case, there are $x,y \in C$ such
that
$$
x+y = \alpha + \delta,
$$
which implies that $\alpha + \delta$ is an effective theta
characteristic. Since $\alpha + \delta$ is even by assumption,
$h^0(C,\alpha + \delta) = 2$, which implies that $C$ is hyperelliptic
and that $\alpha + \delta$ is the hyperelliptic series. \endproof
\medskip

Next we use these cycles to construct various variations of mixed
Hodge
structure whose period maps give lifts
of the period map $\M_3(l) \to \A_3(l)$ to $\calQ\times_{\A_3(l)}
\calP$,
and  to $\calB$ generically.

For each $D \in \Pic^0 C$, one has the extension of mixed Hodge
structure
$$
0 \to H_3(\Pic^1 C;\Z(-1)) \to
H_3(\Pic^1 C, Z_{\alpha,D};\Z(-1)) \to \Z \to 0
$$
where the generator 1 of $\Z$ is the image of any relative class
$[\Gamma]$ with $\partial \Gamma = Z_{\alpha,D}$.  Pushing this
extension out along the projection
$$
H_3(\Pic^1 C) \to QH_3(\Pic^1 C)
$$
one obtains an extension
\begin{equation}\label{Q}
0 \to QH_3(\Pic^1 C;\Z(-1)) \to E \to \Z \to 0.
\end{equation}
This extension determines the same point in
$J(QH_3(\Pic^1 C))\approx Q_1(\Jac C)$
as the cycle $C_x - C_x^-$, and is independent of
the choice of $D$.

Dually, one can consider the extension
\begin{equation}\label{extension}
 0 \to \Z(1) \to H_3(\Pic^1 C - Z_{\alpha,D}; \Z(-1)) \to
H_3(\Pic^1 C; \Z(-1)) \to 0
\end{equation}
which comes from the Gysin sequence. Here the canonical generator
of $\Z(1)$ is the boundary of any 4 ball which is transverse to
$Z_{\alpha,D}$ and has intersection number 1 with it.

\begin{proposition}\label{map}
If $D$ is a point of order 2 in $\Jac C$, then there is a morphism of
Hodge structures
$$
H_1(C,\Z) \to H_3(\Pic^1 C - Z_{\alpha,D}; \Z(-1))
$$
whose composition with the natural map
$$
H_3(\Pic^1 C - Z_{\alpha,D}; \Z(-1)) \to H_3(\Pic^1 C ; \Z(-1))
$$
is  the map $\underline{\phantom{x}}\times [C]$ given by Pontrjagin
product with $[C]$.
\end{proposition}

\proof The exact sequence of Hodge structures
$$
0\to H_1(C;\Z) \stackrel{\times[C]}{\longrightarrow}
 H_3(\Pic^1 C;\Z(-1)) \to QH_3(\Pic^1 C;\Z(-1)) \to 0
$$
induces an exact sequence of Ext groups
$$
0\to \Ext^1(QH_3(\Pic^1 C,\Z(-1)),\Z(1)) \to
\phantom{xxxxxxxxxxxxxxxxxxxxx}
$$
$$
\phantom{xxxxxxxxxxxxx}
\Ext^1(H_3(\Pic^1 C,\Z(-1)),\Z(1)) \to \Ext^1(H_1(C;\Z),\Z(1)) \to 0.
$$
This sequence may be identified naturally with the sequence
$$
0 \to JPH_3(\Pic^1 C;\Z(-1)) \to JH_3(\Pic^1 C;\Z(-1))
\stackrel{\phi}{\to} \Jac C \to 0,
$$
where $\phi$ is the map induced by the morphism of Hodge structures
$H_3(\Jac C;\Z) \to H_1(C;\Z)$ defined by
$$
x\times y \times z \mapsto (x\cdot y)z + (y\cdot z) x + (z\cdot x) y.
$$
Here $x,y,z$ are elements of $H_1(C;\Z)$, $\times$ denotes the
Pontrjagin
product, and  $(\phantom{x}\cdot\phantom{y})$ denotes the
intersection pairing.  To prove the assertion, it suffices to show
that
the image in $\Jac C$ of the class of the extension (\ref{extension})
vanishes. Lefschetz duality gives an isomorphism of (\ref{extension})
with the  the extension
$$
0 \to H_3(\Pic^1 C;\Z(-1)) \to
H_3(\Pic^1 C,Z_{\alpha,D};\Z(-1)) \to \Z \to 0
$$
It follows directly from \cite[(6.7)]{pulte} that the image of this
extension under the map $\phi$ is $K_C - 2 (\alpha + D) = 0$.
\endproof
\medskip

Taking $D=0$ and dividing out by this copy of $H_1(C;\Z)$, we obtain
an extension
\begin{equation}\label{P}
0 \to \Z(1) \to F \to QH_3(\Pic^1 C; \Z(-1)) \to 0.
\end{equation}

Let $\M_3(l,\alpha,\delta)$ be the moduli space of genus 3 curves
with a level $l$ structure, a distinguished even theta characteristic
$\alpha$, and a distinguished point $\delta \in \Pic^0 C$ of order 2
such that $\alpha + \delta$ is also an even theta characteristic.
Denote the universal jacobian over $\M_3(l,\alpha,\delta)$ by
$\calJ \to \M_3(l,\alpha,\delta)$.
The period maps for the extensions (\ref{Q}), (\ref{P}),
respectively,
define maps
$\calJ \to \calQ$ and $\calJ \to \calP$.  These induce a map $\phi$
into their fibered product over $\A_3(l)$ such that the diagram
$$
\matrix{
\calJ & \mapright{\phi} & \calQ \times_{\A_3(l)} \calP \cr
\downarrow & & \downarrow \cr
\M_3(l,\alpha,\delta) & \to & \A_3(l) \cr }
$$
commutes.  Pulling back the $\C^\ast$ bundle
$\calB \to \calQ\times_{\A_3(l)}\calP$ along $\phi$ gives a $\C^\ast$
bundle
$\calL^\ast \to \calJ$.  Denote the corresponding line bundle by
$\calL$.

Denote the relative difference divisor in $\calJ$ by $\calD$ and the
relative theta divisor associated to $\alpha$ by $\vartheta_\alpha$

\begin{lemma}\label{total_chern}
The Chern class of this line bundle is the divisor
 $\calJ$ is $2\calD - 4\vartheta_\alpha$.
\end{lemma}

\proof We construct a meromorphic section of $\calL$.  A point of
$\calJ$ is a curve $C$ and a point $D$ of $\Jac C$.  According to
(\ref{disjoint}), the cycles $Z_\alpha$ and $Z_{\alpha,D}$ have
disjoint supports when $D \not\in \Delta\cup \Theta_\alpha$. In
this case we can consider the mixed Hodge
structure $H_3(\Pic^1 C - Z_{\alpha},
Z_{\alpha,D};\Z(-1))$.\footnote{This
is everywhere a local system. One has to replace it with another
group
when $C$ is hyperelliptic and either $\alpha$ or $\alpha + \delta$ is
the hyperelliptic series. For details, see the footnote on page 887
of
\cite{hain:heights}.}
Dividing this biextension out by the image of the composite
$$
H_1(C,\Z) \to H_3(\Pic^1 C - Z_{\alpha}; \Z(-1))
\to H_3(\Pic^1 C - Z_\alpha, Z_{\alpha,D}; \Z(-1))
$$
of the map of (\ref{extension}) with the natural inclusion produces a
biextension $b_{C,\alpha}$  with weight graded quotients canonically
isomorphic to
$$
\Z, \quad QH_3(\Jac C;\Z(-1)),\quad \Z(1);
$$
the generator 1 of $\Z$ corresponding to any $\Gamma$ with
$\partial \Gamma = Z_{\alpha,D}$,
and the canonical generator $2\pi i$ of $\Z(1)$ being the class of
the
boundary of any small 4 ball which is transverse to $Z_{\alpha}$
and having intersection number 1 with it.  This defines a lift
$$
\tilde{\phi} : \calJ - (\calD\cup \vartheta_\alpha) \to \calB
$$
of the map $\phi : \calJ \to \calQ \times_{\A_3(l)}\calP$.  It
therefore defines a nowhere vanishing holomorphic section $s$ of
$\calL \to \calJ$ on the complement of $\calD \cup
\vartheta_\alpha$.  It follows from \cite[(3.4.3)]{hain:heights} that
$s$ extends to a meromorphic section of $\calL$  on all of $\calJ$.
Consequently, the Chern class of this bundle is supported on the
divisor $\calD \cup \vartheta_\alpha$.  Since the divisors $\calD$
and $\vartheta_\alpha$ are irreducible, the Chern class can be
computed by restricting to a general enough fiber.

With the aid of (\ref{map}), and the formula
\cite[(3.2.11)]{hain:heights}, one can easily show that the height of
the
biextension $b_{C,\alpha}$ equals that of
$$
H_3(\Pic^1 C - Z_\alpha, \Z_{\alpha,D};\Z(-1)).
$$
It follows from
the main theorem of \cite{hain:heights} that the divisor of $s$
restricted to $\Jac C$ is  $2D - 4\Theta_\alpha$ for all $C$. The
result
follows. \endproof
\medskip

The point $\delta$ of order 2 is a section of the bundle $\calJ \to
\M_3(l,\alpha,\delta)$. A lift
$\zeta : \M_3(l,\alpha,\delta) \to \calQ\times_{\A_3(l)}\calP$ of the
period map can be defined by composing $\delta$ with $\phi$.
The pullback of the line bundle $\calL$ along $\delta$ equals the
pullback of the biextension line bundle to $\M_3$.  It follows
that the Chern class of this line bundle is
$2 \delta^\ast (\calD - 2 \vartheta_\alpha).$

\begin{proposition}\label{chern_class}
If $\alpha$ and $\alpha +\delta$ are even theta characteristics,
then the pushforward of the divisor
$\delta^\ast (\calD - 2 \vartheta_\alpha)$ in $\M_3(l,\alpha,\delta)$
to $\M_3(l)$ is $28.35$ times the hyperelliptic locus. Consequently,
the line bundle
$$
\delta^\ast\calL\in \Pic\, \M_3(l,\alpha,\delta)\otimes\Q
$$
is non-trivial.
\end{proposition}

In the proof of this result, we will need the following fact.

\begin{lemma}\label{transverse}
The section $\delta: \M_3(l,\alpha,\delta) \to \calJ$ is transverse
to
the divisors $\calD$ and $\vartheta_\alpha$.
\end{lemma}

\proof  We first prove that $\delta$ intersects $\vartheta_\alpha$
transversally.  We view $\M_3(l,\alpha,\delta)$ as a subvariety of
$\calJ$ via the section $\delta$.  Let $(C,\alpha,\delta)$ be a point
in $\vartheta_\alpha$.  Then, by (\ref{hyperelliptic}), $C$ is
hyperelliptic, and $\alpha + \delta$ is the hyperelliptic series.
So $h^0(\alpha) = 0$ and $h^0(\alpha+\delta) = 2$.  Let $Z_0$ be the
period matrix of $C$ with respect to some symplectic basis of
$H_1(C;\Z)$, and $\theta_\alpha(u,Z)$ the theta function which
defines $\vartheta_\alpha$.  The point $\delta$ of order 2 may be
viewed as a function $\delta(Z)$.  Set $\delta_0 = \delta(Z_0)$.
Since
$\delta_0 \in \Theta_\alpha$, $\theta_\alpha(\delta_0,Z_0)=0$.   We
have to show that there exist $a,b$
such that
$$
{\partial \over \partial Z_{ab}} \theta_\alpha(\delta(Z),Z)|_{Z_0}
\neq 0.
$$
By Riemann's Theorem \cite[p.\ 348]{griffiths-harris}, the
multiplicity of $\delta$ on $\Theta_\alpha$ is
$h^0(\alpha +\delta) = 2$. That is,
\begin{equation}\label{vanishing}
{\partial \theta_\alpha \over \partial u_a} (\delta_0,Z_0) = 0
\end{equation}
for all indices $a$, but there exist indices $a,b$ such that
$$
{\partial^2 \theta_\alpha \over \partial u_a\partial u_b}
(\delta_0,Z_0) \neq 0.
$$
Substituting (\ref{vanishing}) into the chain rule, we have
$$
{\partial \theta_\alpha \over \partial Z_{ab}}(\delta(Z),Z)|_{Z_0}
= \sum_{j=1}^g
{\partial \theta_\alpha\over \partial u_j} (\delta_0,Z_0)
{\partial u_j\over \partial Z_{ab}}(\delta_0,Z_0) +
{\partial \theta_\alpha \over \partial Z_{ab}}(\delta_0,Z_0)
= {\partial \theta_\alpha \over \partial Z_{ab}}(\delta_0,Z_0).
$$
Plugging this into the heat equation, we
obtain the desired result:
$$
{\partial \theta_\alpha \over \partial Z_{ab}}(\delta(Z),Z)|_{Z_0}=
{\partial \theta_\alpha \over \partial Z_{ab}}(\delta_0,Z_0)
= 2\pi i (1 + \delta_{ab})
{\partial^2 \theta_\alpha \over \partial u_a \partial u_b}
(\delta_0,Z_0) \neq 0.
$$

To prove that $\delta$ is transverse to $\calD$, we use an argument
suggested
to us by Nick Katz.  Consider a family of curves
$C \to \spec \C[\epsilon]/(\epsilon^2)$
over the dual numbers.  There is a relative difference divisor
$\Delta\to \spec \C[\epsilon]/(\epsilon^2)$ contained in the Picard
scheme $\Pic^0 C \to \spec \C[\epsilon]/(\epsilon^2)$.  Suppose that
we have a point of order 2
$$
\delta : \spec \C[\epsilon]/(\epsilon^2) \to \Pic^0 C,
$$
defined over the dual numbers  which lies in $\Delta$.  To prove
transversality, it suffices  to show that $C$ is hyperelliptic over
the
dual numbers.  But this is immediate as  $\delta$ gives a 2:1 map
$C \to \P^1$ defined over the dual numbers.  \endproof
\medskip

\noindent{\bf Proof of (\ref{chern_class}).}   Denote the
hyperelliptic
series of a hyperelliptic curve by $H$. It follows from
(\ref{total_chern}) and (\ref{transverse}) that
$$
\delta^\ast \calD = \calH_\Delta\quad\hbox{and}\quad
\delta^\ast \vartheta_\alpha = \calH_0
$$
where
$$
\calH_\Delta = \left\{ (C,\alpha,\delta) : C \hbox{ is hyperelliptic
},
\delta \in \Delta\hbox{ and } h^0(C,\alpha + \delta)
\hbox{ is even}\right\}
$$
and
$$
\calH_0 = \left\{ (C,\alpha,\delta) : C \hbox{ is hyperelliptic },
\alpha \neq H, \alpha + \delta = H \right\}.
$$
If $C$ is hyperelliptic and if $\delta = x - y \in \Delta$, then $x$
and $y$ are distinct Weierstrass points, and
$$
H + \delta = 2y + x -y = x + y
$$
which is an odd theta characteristic. It follows that
$$
\calH_\Delta = \left\{ (C,\alpha,\delta) : C \hbox{ is hyperelliptic
}, \alpha
\neq H, \delta \in \Delta\hbox{ and } h^0(C,\alpha + \delta) \hbox{
even}\right\}.
$$

Apart from the hyperelliptic series, every even theta characteristic
on a
hyperelliptic curve $C$ is of the form
$$
-H + p_1 + p_2 + p_3 + p_4 = - H + q_1 + q_2 + q_3 + q_4
$$
where $p_1,\ldots,p_4, q_1,\ldots,q_4$ are the Weierstrass points.
If $\alpha = -H + p_1 + p_2 + p_3 + p_4 $ and $\alpha + x-y$ is
an even theta characteristic, then it is not difficult to show that
$x-y = q_i - p_j$ for some $i,j$.  So, for each even theta
characteristic
$\alpha \neq H$, there are 16 points $\delta\in \Delta$ such that
$\alpha + \delta$ is also an even theta characteristic. Since there
are
35 even theta characteristics $\alpha \neq H$, it follows that
$\calH_\Delta$ has degree 16.35 over the hyperelliptic locus $\calH$
of $\M_3(l)$.  Since there is only one point $\delta$ of order 2 such
that $\alpha + \delta = H$, $\calH_0$ has degree 35 over $\calH$.

Putting this together we see that the pushforward of
$\delta^\ast\calL$ to $\M_3$
is
$$
\pi_\ast (2\calH_\Delta - 4\calH_0) =
(2.16.35 - 4.35)\calH = 28.35 \calH. \quad \square
$$

We are now ready to prove (\ref{main_3}).  Denote the subgroup of
$\Gamma_3$ which corresponds to $\M_3(l,\alpha,\delta)$ by
$\Gamma_3(l,\alpha,\delta)$, and its intersection with $T_3$
by $T_3(\alpha,\delta)$.
\medskip

\noindent{\bf Proof of (\ref{main_3}).}
Let $N$ be the line bundle over $\A_3(l)$ which is the determinant
of $R^1f_\ast \calO$, where $f: \calJ \to \A_3(l)$ is the universal
abelian variety.  The restriction of this to $\M_3(l)$ is $\calH/9$
\cite[p.~134]{harris_j}.  We shall also denote its pullback to
$\calQ\times_{\A_3(l)}\calP$ by $N$.  It follows from
(\ref{chern_class})
that the line bundle $\calL \otimes N^{\otimes(-9.28.35)}$ pulls back
to
the trivial line bundle over $\M_3(l,\alpha,\delta)$.  There is
therefore
a lift of the period map
$$
\M_3(l,\alpha,\delta) \to
\left(\calL \otimes N^{\otimes(-9.28.35)}\right)^\ast
$$
to the $\C^\ast$ bundle associated to
$\calL \otimes N^{\otimes(-9.28.35)}$.
This induces a group homomorphism
$$
\Gamma_3(l,\alpha,\delta) \to
\pi_1(\left(\calL \otimes N^{\otimes(-9.28.35)}\right)^\ast,\ast).
$$
This last group is an extension
$$
1 \to G_\Z \to
\pi_1(\left(\calL \otimes N^{\otimes(-9.28.35)}\right)^\ast,\ast)
\to L(l) \to 1.
$$
It follows from (\ref{double}) and the fact that $T_3(\alpha,\delta)$
has finite index in $T_3$ that the image of the map
$$
H_1(T_3(\alpha,\delta);\Q) \to H_1(G_\Z;\Q) =
QH_3(\Jac C;\Q) \oplus PH_3(\Jac C;\Q) \approx V(\lambda_3)^2
$$
is the diagonal copy of $V(\lambda_3)$.  The restriction of the
extension
$$
1 \to \Q \to G_\Q\to  V(\lambda^3)^2 \to 1
$$
to the the diagonal is the extension given by the polarization of
$V(\lambda_3)$.  It follows that the homomorphism
$$
\Q \hookrightarrow \Lambda^2 H_1(T_3;\Q) \to \Q
$$
given by evaluating the bracket on the polarization is an
isomorphism, as
claimed.  Second, the element of $H^2(L(l);\Q)$
which corresponds to the extension
$$
0 \to \Q \to H \to L(l) \to 1
$$
constructed from the Torelli group is just the Chern class of the
pullback of
the line bundle $\calL \otimes N^{\otimes(-9.28.35)}$
pulled back to $\A_3(l)$ along the zero section of $\calQ
\times_{\A_3(l)}\calP
\to \A_3(l)$.  This is just $-9.28.35c_1(N)$,
which is nonzero in $H^2(\A_3(l);\Q)$.  Consequently, the extension
is non-trivial as claimed. \endproof
\medskip

\bibliographystyle{plain}

\end{document}